\documentclass[12pt, reqno]{amsart}
\allowdisplaybreaks
\begin{document}

\title[\hfilneg \hfil   inequalities on time scales]
{Estimates of some applicable inequalities on time scales }

\author[D. B. Pachpatte\hfil \hfilneg]
{Deepak B. Pachpatte}

\address{Deepak B. Pachpatte \newline
 Department of Mathematics,
 Dr. B.A.M. University, Aurangabad,
 Maharashtra 431004, India}
\email{pachpatte@gmail.com}


\begin{abstract} 
  The main objective of the paper is to establish explicit estimates on some applicable inequalities in two variables on time scales which can be used in the study of certain qualitative properties of dynamical equations on time scales.
\end{abstract} 

\maketitle

\section{Introduction}
     Many physical, chemical and biological phenomena can be modeled using dynamic equations and study of such problems has enormous potential. In 1988  Stefan Hilger \cite{HIG} in his Ph.D thesis introduced the calculus on time scales which unifies the continuous and discrete analysis. As a response to the diverse need of the applications recently in last decade many authors have studied the properties of solutions of dynamic equations on time scales \cite{And1, And2, And3,ams1, Boh3, Gra1, Hof1, Dbp1, Dbp2, Dbp4, Dbp5, sla1,sla2, wen, mut}. Motivated by the above results in this paper we find  inequalities with explicit estimates which can found to be important tool in the study of dynamical systems on time scales.
     Let  $\mathbb{R}$ denotes the set of real numbers and $\mathbb{T}$ denotes an arbitrary time scale.

          More basic information about time scales calculus can be found in monographs  \cite{Boh1, Boh2}.
Now following \cite{wen,mut}  we give some basic definitions about calculus on time scales in two variables.

 We say that $f:\mathbb{T} \to \mathbb{R}$ is rd-continuous provided $f$ is continuous at each right-dense point of $\mathbb{T}$ and has a finite left sided limit at each left dense point of $\mathbb{T}$. $C_{rd}$ denotes the set of rd-continuous function defined on $\mathbb{T}$. Let $\mathbb{T}_1$ and $\mathbb{T}_2$ be two time scales with at least two points and consider the time scales intervals $\overline {\mathbb{T}}_1  = \left[ {x_0 ,\infty } \right) \cap \mathbb{T}_1$ and $\overline{\mathbb{T}}_2 = \left[ {y_0 ,\infty } \right) \cap \mathbb{T}_2$ for $x_0 \in \mathbb{T}_1 $ and $y_0 \in \mathbb{T}_2 $ and $\Omega=\mathbb{T}_1 \times \mathbb{T}_2$.
 Let $\sigma_1,\rho_1,\Delta_1$ and $\sigma_2,\rho_2,\Delta_2$ denote the forward jump operators, backward jump operators and the delta differentiation operator respectively on $\mathbb{T}_1$
 and $\mathbb{T}_2$. Let $a<b$ be points in $\mathbb{T}_1$, $c<d$ are point in  $\mathbb{T}_2$, $[a,b)$ is the half closed bounded interval in $\mathbb{T}_1$, and $[c,d)$ is the half closed bounded interval in $\mathbb{T}_2$.

 We say that a real valued function $f$ on $\mathbb{T}_1\times\mathbb{T}_2$ at $(t_1,t_2) \in \overline{\mathbb{T}}_1 \times \overline{\mathbb{T}}_2 $ has a $\Delta_1$ partial derivative  $f^{\Delta_1}(t_1,t_2)$ with respect to $t_1$ if for each $\epsilon >0$ there exists a neighborhood $U_{t_1 }$ of $t_1$ such that
\[
\left| {f\left( {\sigma _1 \left( {t_1 } \right),t_2 } \right) - f\left( {s,t_2 } \right) - f^{\Delta _1 } \left( {t_1 ,t_2 } \right)\left( {\sigma _1 \left( {t_1 } \right) - s} \right)} \right| \le \varepsilon \left| {\sigma _1 \left( {t_1 } \right) - s} \right|,
\]
for each $s \in U_{t_1}$,$t_2 \in \mathbb{T}_2$.
We say that  $f$ on $\mathbb{T}_1\times\mathbb{T}_2$ at $(t_1,t_2) \in \overline{\mathbb{T}}_1 \times \overline{\mathbb{T}}_2 $ has a $\Delta_2$ partial derivative  $f^{\Delta_2}(t_1,t_2)$ with respect to $t_2$ if for each $\eta >0$ there exists a neighborhood $U_{t_2 }$ of $t_2$ such that
\[
\left| {f\left( {t_1 ,\sigma _2 \left( {t_2 } \right)} \right) - f\left( {t_1 ,l} \right) - f^{\Delta _2 } \left( {t_1 ,t_2 } \right)\left( {\sigma _2 \left( {t_2 } \right) - l} \right)} \right| \le \eta \left| {\sigma _2 \left( {t_2 } \right) - l} \right|,
\]
for all $l \in U_{t_2}$,$t_1 \in \mathbb{T}_1$.
 The function $f$ is called rd-continuous in $t_2$ if for every $\alpha_1 \in \mathbb{T}_1$, the function $f(\alpha_1,.)$ is rd-continuous on $\mathbb{T}_2$. The function $f$ is called rd-continuous in $t_1$ if for every $\alpha_2 \in \mathbb{T}_2$ the function $f(.,\alpha_2)$ is rd-continuous on $\mathbb{T}_1$.

      The partial delta derivative of $z(x,y)$ for $(x,y) \in \Omega$ with respect to x, y and xy is denoted by $z^{\Delta _1 } \left( {x,y} \right)$, $z^{\Delta _2 } \left( {x,y} \right)$, $z^{\Delta _1 \Delta _2 } \left( {x,y} \right) = z^{\Delta _2 \Delta _1 } \left( {x,y} \right)$.

\section{Main Results}
Now we give our main results.

\paragraph{\textbf{Theorem 2.1}} Let $u,p \in C_{rd} \left( {\Omega ,\mathbb{R}_ +  } \right)$ and $k \ge 0$ is constant.
If
\[
u\left( {x,y} \right) \le k + \int\limits_{x_0 }^x {\int\limits_{s_0 }^s {\int\limits_{y_0 }^y {p\left( {\eta ,\tau } \right)} } } u\left( {\eta ,\tau } \right)\Delta \tau \Delta \eta \Delta s,
\tag{2.1}\]
for $\left( {x,y} \right) \in \Omega $, then
\[
u\left( {x,y} \right) \le ke_{\int\limits_{s_0 }^s {\int\limits_{y_0 }^y {p\left( {\eta ,\tau } \right)\Delta \tau \Delta \eta .} } } \left( {x,x_0 } \right).
\tag{2.2}\]
\paragraph{\textbf{Proof}}  Assume $k>0$. Define a function $w(x,y)$ by right hand side of $(2.1)$, $w(x,0)=w(0,y)=k$, $u(x,y) \le w(x,y)$.
\[
w^{\Delta _2 } \left( {x,y} \right) = \int\limits_{x_0 }^x {\int\limits_{s_0 }^s {p\left( {\eta ,\tau } \right)} } u\left( {\eta ,\tau } \right)\Delta \eta \Delta \tau,
\tag{2.3}\]
\[
w^{\Delta _1 } \left( {x,y} \right) = \int\limits_{x_0 }^x {\int\limits_{y_0 }^y {p\left( {\eta ,\tau } \right)} } u\left( {\eta ,\tau } \right)  \Delta \tau \Delta \eta,
\tag{2.4}\]
\[
w^{\Delta _1 \Delta _1 } \left( {x,y} \right) = \int\limits_{y_0 }^y {p\left( {x,\tau } \right)} u\left( {x,\tau } \right)\Delta \tau,
\tag{2.5}\]
and
\[
w^{\Delta _1 \Delta _1 \Delta _2 } \left( {x,y} \right) = p\left( {x,y} \right)u\left( {x,y} \right) \le p\left( {x,y} \right)w\left( {x,y} \right).
\tag{2.6}\]
From $(2.6)$ and from the facts that $w^{\Delta _1 \Delta _1 }w(x,y) \ge 0, w^{\Delta _1 }w(x,y) \ge 0,w(x,y)>0$ we have
 \[
\frac{{w^{\Delta _1 \Delta _1 \Delta _2 } \left( {x,y} \right)}}{{w\left( {x,y} \right)}} \le p\left( {x,y} \right) + \left[ {\frac{{w^{\Delta _1 \Delta _1 } \left( {x,y} \right)w^{\Delta _2 } \left( {x,y} \right)}}{{w^2 \left( {x,y} \right)}}} \right]
\]
\[
\frac{{w^{\Delta _1 \Delta _1 \Delta _2 } \left( {x,y} \right)}}{{w\left( {x,y} \right)}} \le p\left( {x,y} \right).
\tag{2.7}\]
By keeping $x$ fixed we set $y=\tau$ and then delta integrating with respect to $\tau$ from $y_0$ to $y$ and $w^{\Delta _1 \Delta _1 } \left( {x,y_0 } \right) = 0$
we get
\[
\frac{{w^{\Delta _1 \Delta _1 } \left( {x,y} \right)}}{{w\left( {x,y} \right)}} \le \int\limits_{y_0 }^y {p\left( {x,\tau } \right)\Delta \tau }.
\tag{2.8}\]
From $(2.8)$ and as we have $w^{\Delta _1 } \left( {x,y} \right) \ge 0$, $w(x,y)>0$ we get
\[
\frac{\partial }{{\Delta _1 }}\left( {\frac{{w^{\Delta _1 } \left( {x,y} \right)}}{{w\left( {x,y} \right)}}} \right) \le \int\limits_{y_0 }^y {p\left( {x,\tau } \right)\Delta \tau }.
\tag{2.9}\]
By taking $y$ fixed in $(2.9)$ set $x=\eta$ integrating $\eta$ with respect to $x_0$ to $x$ and $z^{\Delta _1 } \left( {y_0 ,y} \right) = 0$ we have
\[
\frac{{w^{\Delta _1 } \left( {x,y} \right)}}{{w(x,y)}} \le \int\limits_{x_0 }^x {\int\limits_{y_0 }^y {p\left( {\eta ,\tau } \right)\Delta \tau \Delta \eta } }.
\tag{2.10}\]
From $(2.10)$ we get
\[
w\left( {x,y} \right) \le ke_{\int\limits_{s_0 }^s {\int\limits_{y_0 }^y {p\left( {\eta ,\tau } \right)\Delta \tau \Delta \eta } } } \left( {x,x_0 } \right).
\tag{2.11}\]
Using $(2.11)$ in $u(x,y) \le w(x,y)$ we get the result.
\paragraph{\textbf{Theorem 2.2}}
Let $p,q$ be positive  and rd-continuous and $q$ be non decreasing. If
\[
u(x,y) \le q\left( {x,y} \right) + \int\limits_{x_0 }^x {\int\limits_{s_0 }^s {\int\limits_{y_0 }^y {p\left( {\eta ,\tau } \right)u\left( {\eta ,\tau } \right)\Delta \tau \Delta \eta \Delta s} } },
\tag{2.12}\]
for $\left( {x,y} \right) \in \Omega $, then
\[
u(x,y) \le q\left( {x,y} \right)e_{\int\limits_{s_0 }^s {\int\limits_{y_0 }^y {p\left( {\eta ,\tau } \right)u\left( {\eta ,\tau } \right)\Delta \tau \Delta \eta } } } \left( {x,x_0 } \right).
\tag{2.13}\]
\paragraph{\textbf{Proof} } Let $q(x,y) \ge 0$ for $\left( {x,y} \right) \in \Omega $ . Then from $(2.12)$, it is easy to see that
\[
\frac{{u(x,y)}}{{q\left( {x,y} \right)}} \le 1 + \int\limits_{x_0 }^x {\int\limits_{s_0 }^s {\int\limits_{y_0 }^y {p\left( {\eta ,\tau } \right)\frac{{u\left( {\eta ,\tau } \right)}}{{q\left( {\eta ,\tau } \right)}}\Delta \tau \Delta \eta } } } \Delta s.
\tag{2.14}\]
Now an application of inequality in Theorem 2.1 gives the result $(2.13)$.

\paragraph{\textbf{Theorem 2.3.}}  Let $u,g,h,p \in C_{rd} \left( {\Omega ,\mathbb{R}_ +  } \right)$ and $L \in C_{rd} \left( {\Omega  \times \mathbb{R}_ +  ,\mathbb{R}_ +  } \right)$
\[
 0\le L(x,y,u)-L(x,y,v) \le H(x,y,v)(u-v),
 \tag{2.15}\]
 and $u\ge v \ge 0$, where $H \in C_{rd} \left( {\Omega  \times \mathbb{R}_ +  ,\mathbb{R}_ +  } \right)$. If

  \[
u\left( {x,y} \right) \le g\left( {x,y} \right) + h\left( {x,y} \right)\int\limits_{x_0 }^x {\int\limits_{s_0 }^s {\int\limits_{y_0 }^y {L\left( {\eta ,\tau ,u\left( {\eta ,\tau } \right)} \right)} } } \Delta \tau \Delta \eta \Delta s,
\tag{2.16}\]
for $(x,y) \in \Omega $ then
\begin{align*}
u\left( {x,y} \right)
& \le g\left( {x,y} \right) + h\left( {x,y} \right)\left( {\int\limits_{x_0 }^x {\int\limits_{s_0 }^s {\int\limits_{y_0 }^y {L\left( {\eta ,\tau ,g\left( {\eta ,\tau } \right)} \right)\Delta \tau \Delta \eta \Delta s} } } } \right) \\
&\times e_{\int\limits_{s_0 }^s {\int\limits_{y_0 }^y {H\left( {\eta ,\tau ,g\left( {\eta ,\tau } \right)} \right)h\left( {\eta ,\tau } \right)\Delta \tau \Delta \eta } } } \left( {x,x_0 } \right),
\tag{2.17}
\end{align*}
for $(x,y) \in \Omega $.

\paragraph{\textbf{Proof}} Define a function $w(x,y)$ by
\[
w\left( {x,y} \right) = \int\limits_{x_0 }^x {\int\limits_{s_0 }^s {\int\limits_{y_0 }^y {L\left( {\eta ,\tau ,u\left( {\eta ,\tau } \right)} \right)\Delta \tau \Delta \eta \Delta s} } },
\tag{2.18}\]
then $w\left( {x,y_0} \right) = w\left( {x_0,y } \right) = 0$ and inequality $(2.16)$ becomes
\begin{align*}
w\left( {x,y} \right)
&\le \int\limits_{x_0 }^x {\int\limits_{s_0 }^s {\int\limits_{y_0 }^y {\left\{ {L\left( {\eta ,\tau ,g\left( {\eta ,\tau } \right) + h\left( {\eta ,\tau } \right)w\left( {\eta ,\tau } \right)} \right)} \right.} } }  \\
&\left. { - L\left( {\eta ,\tau ,g\left( {\eta ,\tau } \right)} \right) + L\left( {\eta ,\tau ,g\left( {\eta ,\tau } \right)} \right)} \right\}\Delta \tau \Delta \eta \Delta s \\
&\le \int\limits_{x_0 }^x {\int\limits_{s_0 }^s {\int\limits_{y_0 }^y {H\left( {\eta ,\tau ,g\left( {\eta ,\tau } \right)} \right)h\left( {\eta ,\tau } \right)w\left( {\eta ,\tau } \right)\Delta \tau \Delta \eta \Delta s} } } \\
& + \int\limits_{x_0 }^x {\int\limits_{s_0 }^s {\int\limits_{y_0 }^y {L(\eta ,\tau ,g\left( {\eta ,\tau } \right))\Delta \tau \Delta \eta \Delta s} } }. 
\tag{2.19}
\end{align*}

It can be easily seen that the first term on the right hand side of $(2.19)$ is nonnegative and non decreasing. Now suitable application of Theorem $2.2$ to $(2.19)$  we get $(2.17)$.

\paragraph{\textbf{Theorem 2.4}} Let $u,g,h,p$ be as in theorem 2.3.
If
\[
u\left( {x,y} \right) \le g\left( {x,y} \right) + h\left( {x,y} \right)\int\limits_{x_0 }^x {\int\limits_{s_0 }^s {\int\limits_{y_0 }^y {p\left( {\eta ,\tau } \right)} } } u\left( {\eta ,\tau } \right)\Delta \tau \Delta \eta \Delta s,
\tag{2.19}\]
for $(x,y) \in \Omega $, then
\begin{align*}
u\left( {x,y} \right)
&\le g\left( {x,y} \right) + h\left( {x,y} \right)\left( {\int\limits_{x_0 }^x {\int\limits_{s_0 }^s {\int\limits_{y_0 }^y {p\left( {\eta ,\tau } \right)g\left( {\eta ,\tau } \right)\Delta \tau \Delta \eta \Delta s} } } } \right) \\
&\times e_{\int\limits_{s_0 }^s {\int\limits_{y_0 }^y {p\left( {\eta ,\tau } \right)h\left( {\eta ,\tau } \right)\Delta \tau \Delta \eta } } } \left( {x,x_0 } \right),
\tag{2.20}
\end{align*}
for $(x,y) \in \Omega $.

\paragraph{\textbf{Proof}} Now putting $L\left( {\eta ,\tau ,u\left( {\eta ,\tau } \right)} \right) = p\left( {\eta ,\tau } \right)u\left( {\eta ,\tau } \right)$ in above  Theorem $2.3$, we get the result.

\section{Applications}
Consider integral equation on time scales of the form
\[
u\left( {x,y} \right) = g\left( {x,y} \right) + \int\limits_{x_0 }^x {\int\limits_{s_0 }^s {\int\limits_{y_0 }^y {K\left( {x,y,\eta ,\tau ,u\left( {\eta ,\tau } \right)} \right)} } } \Delta \tau \Delta \eta \Delta s,
\tag{3.1}\]
where $u$ is unknown function to be found for given $g \in C_{rd} \left( {\Omega ,\mathbb{R}} \right)$ and $K \in C_{rd} \left( {\Omega  \times \Omega  \times \mathbb{R},\mathbb{R}} \right)$.

Now we give the estimates on the solutions of equation $(3.1)$.

\paragraph{\textbf{Theorem 3.1}}   Let $g,K$ in $(3.1)$ satisfy the condition
\[
\left| {K\left( {x,y,\eta ,\tau ,u} \right) - K\left( {x,y,\eta ,\tau ,v} \right)} \right| \le q\left( {x,y} \right)r\left( {\eta ,\tau } \right)\left| {u - v} \right|,
\tag{3.2}\]
\[
\left| {g\left( {x,y} \right) + \int\limits_{x_0 }^x {\int\limits_{s_0 }^s {\int\limits_{y_0 }^y {K\left( {x,y,\eta ,\tau ,0} \right)} } } \Delta \tau \Delta \eta \Delta s} \right| \le p\left( {x,y} \right),
\tag{3.3}\]
where $p,q,r \in C_{rd} \left( {\Omega  ,\mathbb{R}_ +  } \right)$. If $u(x,y)$ is solution of $(3.1)$ for $\left( {x,y} \right) \in \Omega  $, then
\begin{align*}
\left| {u\left( {x,y} \right)} \right|
&\le p\left( {x,y} \right) + q\left( {x,y} \right)\left( {\int\limits_{x_0 }^x {\int\limits_{s_0 }^s {\int\limits_{y_0 }^y {r\left( {\eta ,\tau } \right)p\left( {\eta ,\tau } \right)} } } \Delta \tau \Delta \eta \Delta s} \right) \\
& \times e_{\int\limits_{s_0 }^s {\int\limits_{y_0 }^y {r\left( {\eta ,\tau } \right)q\left( {\eta ,\tau } \right)\Delta \tau \Delta \eta } } } \left( {x,x_0 } \right),
\tag{3.4}
\end{align*}
for $(x,y) \in \Omega $.

\paragraph{\textbf{Proof}} We have $u(x,y)$ as solution of $(3.1)$ for $(x,y) \in \Omega$. We have
\begin{align*}
 u\left( {x,y} \right)
 &\le \left| {g\left( {x,y} \right) + \int\limits_{x_0 }^x {\int\limits_{s_0 }^s {\int\limits_{y_0 }^y {K\left( {x,y,\eta ,\tau ,0} \right)} \Delta \tau \Delta \eta \Delta s} } } \right| \\
 &+ \int\limits_{x_0 }^x {\int\limits_{s_0 }^s {\int\limits_{y_0 }^y {\left| {K\left( {x,y,\eta ,\tau ,u\left( {\eta ,\tau } \right)} \right) - K\left( {x,y,\eta ,\tau ,0} \right)} \right|} \Delta \tau \Delta \eta \Delta s} }  \\
 &\le p\left( {x,t} \right) + q\left( {x,t} \right)\int\limits_{x_0 }^x {\int\limits_{s_0 }^s {\int\limits_{y_0 }^y {r\left( {\eta ,\tau } \right)\left| {u\left( {\eta ,\tau } \right)} \right|} \Delta \tau \Delta \eta \Delta s} }.
\tag{3.5}
\end{align*}

Now applying Theorem 2.4 to $(3.5)$ gives $(3.4)$.\\
\\

A function $u \in C_{rd} \left( {\Omega ,\mathbb{R}} \right)$  is called $\epsilon$ approximate solution of equation $(3.1)$ if their exists a $\epsilon \ge 0$ such that
\[
\left| {u\left( {x,y} \right) - \left\{ {g\left( {x,y} \right) + \int\limits_{x_0 }^x {\int\limits_{s_0 }^s {\int\limits_{y_0 }^y {K\left( {x,y,\eta ,\tau ,u\left( {\eta ,\tau } \right)} \right)\Delta \tau \Delta \eta \Delta s} } } } \right\}} \right| \le \epsilon,
\tag{3.6}\]
for $(x,y) \in \Omega $.

Now we estimate the difference between two approximate solution of $(3.1)$.

\paragraph{\textbf{Theorem 3.2}}
Let $u_i({x,y}) (i=1,2)$ be  $\epsilon_i$ approximate solutions of $(3.1)$ for $(x,y) \in \Omega $. Suppose function $K$ satisfies the condition $(3.2)$. Then
\begin{align*}
 \left| {u_1 \left( {x,y} \right) - u_2 \left( {x,y} \right)} \right|
 &\le \left( {\varepsilon _1  + \varepsilon _2 } \right)\left[ {1 + } \right.q\left( {x,y} \right)\left( {\int\limits_{x_0 }^x {\int\limits_{s_0 }^s {\int\limits_{y_0 }^y {r\left( {\eta ,\tau } \right)w\left( {\eta ,\tau } \right)\Delta \tau \Delta \eta \Delta s} } } } \right) \\
 &\times e_{\int\limits_{s_0 }^s {\int\limits_{y_0 }^y {r\left( {\eta ,\tau } \right)q\left( {\eta ,\tau } \right)\Delta \tau \Delta \eta } } } \left( {x,x_0 } \right),
\tag{3.7}
\end{align*}
for $(x,y) \in \Omega $.

\paragraph{\textbf{Proof}}
Since $u_i({x,y}) (i=1,2)$ be $\epsilon_i$ approximate solutions of $(3.1)$ we get
\[
\left| {u_i \left( {x,y} \right) - \left\{ {g\left( {x,y} \right) + \int\limits_{x_0 }^x {\int\limits_{s_0 }^s {\int\limits_{y_0 }^y {K\left( {x,y,\eta ,\tau ,u_i \left( {\eta ,\tau } \right)} \right)\Delta \tau \Delta \eta \Delta s} } } } \right\}} \right| \le \varepsilon _i.
\tag{3.8}\]
From $(3.8)$ and using the inequalities
\[
\left| {v - \overline v } \right| \le \left| v \right| + \left| {\overline v } \right|,\,\,\,\left| v \right| - \left| {\overline v } \right|\, \le \left| {v - \overline v } \right|,
\tag{3.9}\]
we have

\begin{align*}
\epsilon _1  + \epsilon _2
&\ge \left| {u_1 \left( {x,y} \right) - \left\{ {g\left( {x,y} \right) + \int\limits_{x_0 }^x {\int\limits_{s_0 }^s {\int\limits_{y_0 }^y {K\left( {x,y,\eta ,\tau ,u_1 \left( {\eta ,\tau } \right)} \right)\Delta \tau \Delta \eta \Delta s} } } } \right\}} \right| \\
&+ \left| {u_2 \left( {x,y} \right) - \left\{ {g\left( {x,y} \right) + \int\limits_{x_0 }^x {\int\limits_{s_0 }^s {\int\limits_{y_0 }^y {K\left( {x,y,\eta ,\tau ,u_2 \left( {\eta ,\tau } \right)} \right)\Delta \tau \Delta \eta \Delta s} } } } \right\}} \right| \\
&\ge \left| {\left[ {u_1 \left( {x,y} \right) - \left\{ {g\left( {x,y} \right) + \int\limits_{x_0 }^x {\int\limits_{s_0 }^s {\int\limits_{y_0 }^y {K\left( {x,y,\eta ,\tau ,u_1 \left( {\eta ,\tau } \right)} \right)\Delta \tau \Delta \eta \Delta s} } } } \right\}} \right]} \right. \\
& - \left. {\left[ {u_2 \left( {x,y} \right) - \left\{ {g\left( {x,y} \right) + \int\limits_{x_0 }^x {\int\limits_{s_0 }^s {\int\limits_{y_0 }^y {K\left( {x,y,\eta ,\tau ,u_2 \left( {\eta ,\tau } \right)} \right)\Delta \tau \Delta \eta \Delta s} } } } \right\}} \right]} \right| \\
&\ge \left| {u_1 \left( {x,y} \right) - u_2 \left( {x,y} \right)} \right| - \left| {\int\limits_{x_0 }^x {\int\limits_{s_0 }^s {\int\limits_{y_0 }^y {K\left( {x,y,\eta ,\tau ,u_1 \left( {\eta ,\tau } \right)} \right)\Delta \tau \Delta \eta \Delta s} } } } \right. \\
&\left. { - \int\limits_{x_0 }^x {\int\limits_{s_0 }^s {\int\limits_{y_0 }^y {K\left( {x,y,\eta ,\tau ,u_2 \left( {\eta ,\tau } \right)} \right)\Delta \tau \Delta \eta \Delta s} } } } \right|. \\
\tag{3.10}
\end{align*}
Let $w\left( {x,y} \right) = \left| {u_1 \left( {x,y} \right) - u_2 \left( {x,y} \right)} \right|$ for any $(x,y) \in \Omega $.
From $(3.10)$ and using $(3.2)$, we have
\begin{align*}
w\left( {x,y} \right)
&\le \left( {\epsilon _1  + \epsilon _2 } \right) \\
&+ \int\limits_{x_0 }^x {\int\limits_{s_0 }^s {\int\limits_{y_0 }^y {\left| {K\left( {x,y,\eta ,\tau ,u_1 \left( {\eta ,\tau } \right)} \right) - K\left( {x,y,\eta ,\tau ,u_2 \left( {\eta ,\tau } \right)} \right)} \right|} } } \Delta \tau \Delta \eta \Delta s \\
&\le \left( {\varepsilon _1  + \varepsilon _2 } \right) + q\left( {x,y} \right)\int\limits_{x_0 }^x {\int\limits_{s_0 }^s {\int\limits_{y_0 }^y {r\left( {\eta ,\tau } \right)w\left( {\eta ,\tau } \right)\Delta \tau \Delta \eta \Delta s} } } . 
\tag{3.11}
\end{align*}

Now an using inequality in Theorem $2.3$ yields the result.

\paragraph{\textbf{Remark}}
 In case $u_1 \left( {x,y} \right)$ is a solution of $(3.1)$ then we have $\epsilon_1 = 0$ and from $(3.7)$ we have $u_2 ({x,y})\rightarrow u_1 ({x,y})$ as $\epsilon_2 \rightarrow 0$. If we put $\epsilon_1=\epsilon_2=0$ in $(3.7)$ then the uniqueness of solution of equation $(3.1)$ is established.

{\bf ACKNOWLEDGEMENTS.} The research in the present paper is supported by Science and Engineering Research Board(SERB, New Delhi, India),  File No. SR/S4/MS-861/13.

\end{document}